
\def\LaTeX{%
  \let\Begin\begin
  \let\End\end
  \let\salta\relax
  \let\finqui\relax
  \let\futuro\relax}

\def\UK{\def\our{our}\let\sz s}
\def\USA{\def\our{or}\let\sz z}

\UK



\LaTeX

\USA


\salta

\documentclass[twoside,12pt]{article}
\setlength{\textheight}{24cm}
\setlength{\textwidth}{16cm}
\setlength{\oddsidemargin}{2mm}
\setlength{\evensidemargin}{2mm}
\setlength{\topmargin}{-15mm}
\parskip2mm


\usepackage[usenames,dvipsnames]{color}
\usepackage{amsmath}
\usepackage{amsthm}
\usepackage{amssymb}
\usepackage[mathcal]{euscript}

%
%


\definecolor{viola}{rgb}{0.3,0,0.7}
\definecolor{ciclamino}{rgb}{0.5,0,0.5}

\def\gianni #1{{\color{red}#1}}
\def\pier #1{{\color{red}#1}}
\def\juerg #1{{\color{Green}#1}}
\def\revis #1{{\color{blue}#1}}

\def\revis #1{#1}

\def\pier #1{#1}
\def\juerg #1{#1}
\def\gianni #1{#1}




\bibliographystyle{plain}


%

\finqui

\def\Beq{\Begin{equation}}
\def\Eeq{\End{equation}}
\def\Bsist{\Begin{eqnarray}}
\def\Esist{\End{eqnarray}}

\def\Bthm{\Begin{theorem}}
\def\Ethm{\End{theorem}}
\def\Blem{\Begin{lemma}}
\def\Elem{\End{lemma}}

\def\Bdim{\Begin{proof}}
\def\Edim{\End{proof}}
\def\Bcenter{\Begin{center}}
\def\Ecenter{\End{center}}
\let\non\nonumber




\def\step #1 \par{\medskip\noindent{\bf #1.}\quad}


\def\Lip{Lip\-schitz}
\def\Holder{H\"older}
\def\frechet{Fr\'echet}
\def\aand{\quad\hbox{and}\quad}

\def\lhs{left-hand side}
\def\rhs{right-hand side}
\def\sfw{straightforward}


\def\regulariz{regulari\sz}

\def\nbh{neighb\our hood}
\def\bhv{behavi\our}


\def\multibold #1{\def\arg{#1}%
  \ifx\arg\pto \let\next\relax
  \else
  \def\next{\expandafter
    \def\csname #1#1#1\endcsname{{\bf #1}}%
    \multibold}%
  \fi \next}

\def\pto{.}

\def\multical #1{\def\arg{#1}%
  \ifx\arg\pto \let\next\relax
  \else
  \def\next{\expandafter
    \def\csname cal#1\endcsname{{\cal #1}}%
    \multical}%
  \fi \next}


\def\multimathop #1 {\def\arg{#1}%
  \ifx\arg\pto \let\next\relax
  \else
  \def\next{\expandafter
    \def\csname #1\endcsname{\mathop{\rm #1}\nolimits}%
    \multimathop}%
  \fi \next}

\multibold
qwertyuiopasdfghjklzxcvbnmQWERTYUIOPASDFGHJKLZXCVBNM.

\multical
QWERTYUIOPASDFGHJKLZXCVBNM.

\multimathop
ad dist div dom meas sign supp .


\def\accorpa #1#2{\eqref{#1}--\eqref{#2}}
\def\Accorpa #1#2 #3 {\gdef #1{\eqref{#2}--\eqref{#3}}%
  \wlog{}\wlog{\string #1 -> #2 - #3}\wlog{}}


\def\graffe #1{\mathopen\{#1\mathclose\}}

\def\<#1>{\mathopen\langle #1\mathclose\rangle}
\def\norma #1{\mathopen \| #1\mathclose \|}

\def\[#1]{\mathopen\langle\!\langle #1\mathclose\rangle\!\rangle}

\def\iot {\int_0^t}

\def\intQt{\int_{Q_t}}

\def\iO{\int_\Omega}

\def\dt{\partial_t}

\def\cpto{\,\cdot\,}

\def\sfrac #1#2{{\textstyle\frac {#1}{#2}}}

\def\checkmmode #1{\relax\ifmmode\hbox{#1}\else{#1}\fi}
\def\aeO{\checkmmode{a.e.\ in~$\Omega$}}
\def\aeQ{\checkmmode{a.e.\ in~$Q$}}

\def\aaO{\checkmmode{for a.a.~$x\in\Omega$}}

\def\aat{\checkmmode{for a.a.~$t\in(0,T)$}}


\def\erre{{\mathbb{R}}}




\def\genspazio #1#2#3#4#5{#1^{#2}(#5,#4;#3)}
\def\spazio #1#2#3{\genspazio {#1}{#2}{#3}T0}

\def\L {\spazio L}
\def\H {\spazio H}

\def\C #1#2{C^{#1}([0,T];#2)}


\def\Lx #1{L^{#1}(\Omega)}
\def\Hx #1{H^{#1}(\Omega)}

\def\LQ #1{L^{#1}(Q)}

\def\Ldue{\Lx 2}

\def\Huno{\Hx 1}


\def\LQ #1{L^{#1}(Q)}


\let\badtheta\theta
\let\theta\vartheta
\let\eps\varepsilon
\let\phi\varphi

\let\TeXchi\chi                         
\newbox\chibox
\setbox0 \hbox{\mathsurround0pt $\TeXchi$}
\setbox\chibox \hbox{\raise\dp0 \box 0 }
\def\chi{\copy\chibox}



\def\katheta{\kappa(\theta)}
\def\kmin{\kappa_*}
\def\kmax{\kappa^*}
\def\keps{\kappa_\alpha}
\def\Geps{G_\alpha}

\def\thetastar{\theta_*}
\def\thetaz{\theta_0}
\def\thetaeps{\theta_\eps}
\def\thetabar{\overline\theta}
\def\thetan{\theta_n}
\def\thetabarn{\overline\theta_n}
\def\sigman{\sigma_n}
\def\ueps{u_\eps}
\def\sigmaeps{\sigma_\eps}

\def\dK{d_K}
\def\PK{P_{\!K}}
\def\QK{Q_K}
\def\IH{I_{\!H}}
\def\dKeps{d^{\,\eps}_K}
\def\rhostar{\rho^*}
\def\Tn{T_n}
\def\Tnp{T_{n,p}}
\def\PI{p_I}
\def\QI{q_I}
\def\psieps{\psi_\eps}
\def\phieps{\phi_\eps}
\def\Tstar{T^*}
\def\Tg{T_\gamma}
\def\Teps{T^*_{\eps,\alpha}}
\def\tz{t_0}
\def\tu{t_1}
\def\veps{v_\eps}

\def\Vp{V^*}

\def\normaH #1{\norma{#1}_H}

\def\ceps{c_\eps}
\def\ca{c_\alpha}

\Begin{document}


%
\title{\pier{Constrained evolution\\[0.3cm] 
  for a quasilinear parabolic equation}}
\author{}
\date{}
\maketitle
\Bcenter
\vskip-1cm
{\large\sc Pierluigi Colli$^{(1)}$}\\
{\normalsize e-mail: {\tt pierluigi.colli@unipv.it}}\\[.25cm]
{\large\sc Gianni Gilardi$^{(1)}$}\\
{\normalsize e-mail: {\tt gianni.gilardi@unipv.it}}\\[.25cm]
{\large\sc J\"urgen Sprekels$^{(2)}$}\\
{\normalsize e-mail: {\tt sprekels@wias-berlin.de}}\\[.45cm]
$^{(1)}$
{\small Dipartimento di Matematica ``F. Casorati'', Universit\`a di Pavia}\\
{\small via Ferrata 1, 27100 Pavia, Italy}\\[.2cm]
$^{(2)}$
{\small Department of Mathematics}\\
{\small Humboldt-Universit\"at zu Berlin}\\
{\small Unter den Linden 6, 10099 Berlin, Germany}\\[2mm]
{\small and}\\[2mm]
{\small Weierstrass Institute for Applied Analysis and Stochastics}\\
{\small Mohrenstrasse 39, 10117 Berlin, Germany}
\Ecenter
\Begin{abstract}
\pier{In the present contribution\revis{,} a feedback control law is studied for a quasilinear parabolic equation.  First, we prove the well-posedness and some regularity results for the Cauchy--Neumann problem for this equation, modified by adding an extra term which is a multiple of the subdifferential of the distance function from a closed convex set of the space of square-integrable functions. 
Then, we consider convex sets of obstacle or double-obstacle type \gianni{and prove rigorously the following property: if the factor in front of the feedback control is sufficiently large, then the solution reaches the convex set within a finite time and then moves inside it.}}
\End{abstract}
\salta
\pagestyle{myheadings}
\newcommand\testopari{\sc Colli \ --- \ Gilardi \ --- \ Sprekels}
\newcommand\testodispari{\sc \pier{Constrained evolution
  for a quasilinear parabolic equation}}
\markboth{\testodispari}{\testopari}
\finqui
%

\section{Introduction}
\label{Intro}
\setcounter{equation}{0}

{{} 

A notewhorty interest has arisen in the mathematical literature of the last twenty years
for \revis{sliding mode} control (SMC) problems. SMC is considered as a main tool 
for the systematic design of robust controllers for nonlinear complex
dynamical systems operating under uncertainty. The advantage of SMC 
is the separation of the motion of the overall system
in independent partial components with a lower dimension.

The design of feedback control systems with sliding modes is related to the selection 
of suitable control functions enforcing motions along ad-hoc manifolds. 
Hence, a manifold of lower dimension (called the sliding manifold) 
has to be identified such that the original 
system\revis{,} restricted to this sliding manifold\revis{,} has a desired \bhv; then,  
one acts on the system through the control in order to constrain 
the evolution on it, that is\revis{,} to design \revis{an} SMC law that forces
the trajectories of the system to reach the sliding surface and remain on~it
(see, e.g., \cite{PiUs} and references therein). 

Known methods developed for ODEs (cf., e.g., \cite{LO02}) have been recently 
extended to the control of infinite-dimensional dynamical systems. 
For instance, we mention some papers dealing 
with SMC for semilinear PDEs: \cite{CRS11} deals with the stabilization problem 
of a one-dimensional unstable heat conduction system (rod) 
ruled by \revis{a} parabolic partial differential equation 
with a Dirichlet type actuator
from one of the boundaries; in \cite{Levaggi13}\revis{,} \revis{an} SMC law is studied for
a class of parabolic systems where the control acts through a Neumann boundary 
condition; a delay-independent SMC strategy was proposed in \cite{XLGK13} 
to control a class of quasilinear parabolic PDE systems with time-varying delay. 

The recent paper \cite{BCGMR}, in which two of the three authors of this note are involved, 
faces different kinds of SMC problems for a standard phase field system. This system couples  two parabolic equations in \revis{terms} of the variables temperature and order parameter. Sliding manifolds are considered both for a linear combination of variables and just for the order parameter. It is shown that the chosen SMC laws force the system to reach 
within finite time the sliding manifold.
In particular, the control law is nonlocal in space for two of the examined problems. 
When reporting the related results during a conference, the third author of this note 
observed that it was analogously interesting, from the point of view of control problems, 
to force variables to reach not a single elected manifold but instead a closed convex 
\revis{subset} of the space in which the variables \revis{still can} move.

\revis{We} started to think about \revis{it} and, in order to develop this idea, \revis{decided} to
argue first on a single nonlinear equation, of course an evolutionary equation, and 
of parabolic type like
\Beq
  \dt\theta - \div \bigl( \katheta \nabla\theta \bigr) = f
  \quad \hbox{in $Q:=\Omega\times(0,T)$},
   \label{Iparab}
\Eeq
which fits into a well\revis{-}established subject (let us quote some monographs, i.e., 
\cite{Barbu, Brezis, DiB,Lions,Roub, Show, Zhe}). With the aim of discussing 
existence and dynamics of solutions in the framework of the space $L^2(\Omega)$, and being 
interested to reach a closed convex $K\subset L^2(\Omega) $ in \revis{finite} time, 
a feedback control is added to the equation \eqref{Iparab} by considering
\Beq
  \dt\theta - \div \bigl( \katheta \nabla\theta \bigr) + \rho\partial\dK(\theta) \ni f ,
   \label{Imodif}
\Eeq
where $\partial\dK$ is the subdifferential of the distance function $\dK$ 
associated \revis{with}~$K$ and $\rho$ is a positive parameter, to be suitably chosen in order 
to force the solution to enter the convex set (if it is not already \revis{inside}).
\gianni{We complement \eqref{Imodif} by homogeneous Neumann boundary conditions 
and an initial condition like $\theta(0)=\thetaz$.}

\gianni{It is worth noting that our goal is not the mere reaching of the convex set.
We want to allow an evolution inside it, indeed.
On the other hand, it is known that some single elements can be reached in final time
by a controlled evolution ruled by easier feedback control laws.
For instance, by assuming $\kappa$ to be a constant and $f=0$ in~\eqref{Imodif},
if we replace the subdifferential $\partial\dK(\theta)$ 
by $\sign\theta$, where $\sign$ is the usual maximal monotone graph related to the standard $\sign$ function,
and thus write
\Beq
  \dt\theta - \Delta\theta + \rho\sign\theta \ni 0 ,
   \non 
\Eeq
we obtain the closed-loop system (5.29) of \cite[p.~203]{Barbu},
and it is proved there that the trajectory reaches zero in a final time.
In this case the evolution of $\theta$ is completely stopped.
On the contrary, if $K$ is the convex set of the nonnegative functions $v\in\Ldue$, $\rho=1$,
$\thetaz=-1$ and $\theta$ is the space independent function given by $\theta(t)=t-1$ for $t\in[0,T]$,
then, by using the forthcoming formulas \eqref{defQK} and~\eqref{derdK}, 
one can check that \eqref{Imodif} holds with the space independent function $f:(0,T)\to\erre$ defined a.e.~by
\Beq
  f(t) = 1 - |\Omega|^{-1/2}
  \quad \hbox{if $t<1$}
  \aand
  f(t) = 1
  \quad \hbox{if $t\geq1$}
  \non
\Eeq
where $|\Omega|$ is the measure of~$\Omega$.
Thus, $K$~is reached at the time $t=1$ and an evolution continues in $K$ for $t>1$.}

We point out that 
\eqref{Imodif} has the structure of an evolution inclusion 
(cf.~\cite{Barbu, Brezis}) but it 
not a standard variational inequality constraining the solution 
to stay inside the convex set.
On the other hand, one may exert the control on \eqref{Imodif} 
via the parameter $\rho$: 
we can imagine that \revis{the larger} the coefficient $\rho$ is, 
\revis{the faster} the solution will enter the convex. 

In our approach, we can deal with convex sets in $\Lx2$ of obstacle and double-obstacle type.
We are able to treat \revis{these} cases: of course, the analysis is not straightforward, as the reader will see, and it will \revis{also become} 
clear why obstacle convex sets are suitable for us.
This paper is a first attempt to approach a somehow new field of investigation, \revis{and so} we ask the reader to be \revis{generous}: indeed, to the best of our knowledge, at the moment we are not aware of other similar essays. 

\gianni{We} discuss existence and uniqueness of the solution to the 
initial-boundary value problem without any restriction on
$\rho>0 $ and on the \revis{nonempty} closed convex $K$ of $\Lx2$. Then, we focus on 
convex sets of obstacle type and prove that for a sufficiently large $\rho$ 
the solution $\theta $ will reach the convex set in \revis{finite time}.
}

The paper is organized as follows.
In the next section, we list our assumptions, state the problem in a precise form
and present our results.
In Section~\ref{APPROX}, we introduce and solve an approximating problem
which is useful to construct a solution to the problem at hand.
The proofs of our results are then given in Sections~\ref{PROOFS} and~\ref{GOAL}, mainly,
since the Appendix is just devoted to establish a technical lemma.


\section{Statement of the problem and results}
\label{STATEMENT}
\setcounter{equation}{0}

In this section, we describe the problem under study and present our results.
First, we assume $\Omega$ to be a smooth bounded open subset of~$\erre^d$.
Now, we specify the assumptions on the structure of our system.
As for~$\kappa$, we assume that
\Beq
  \kappa : \erre \to \erre
  \quad \hbox{is continuous, nonnegative and bounded,}
  \label{hpkappa}
\Eeq
and set
\Beq
  \kmin := \inf \kappa , \quad
  \kmax := \sup \kappa
  \aand
  G(r) := \int_0^r \kappa(s) \, ds
  \quad \hbox{for $r\in\erre$}.
  \label{defG}
\Eeq
The last condition defines the function $G:\erre\to\erre$, and we suppose that
\Beq
  \hbox{$G$ is strictly increasing}.
  \label{hpG}
\Eeq
This assumption is satisfied if and only if the set where the nonnegative function $\kappa$ vanishes
has an empty interior, and it clearly holds if $\kmin>0$.
In \juerg{this} case, the inverse function
$G^{-1}$ is \Lip\ continuous and not only continuous
on its domain $(\inf G,\sup G)$.
Next, even for a future convenience, we introduce the spaces
\Beq
  H := \Ldue 
  \aand
  V := \Huno, 
  \label{defspazi}
\Eeq
and we endow them with their standard norms.
It is understood that $H$ is embedded in~$\Vp$, the dual space of~$V$,
in the usual way, i.e., \juerg{such} that
$\<u,v>=\iO uv\,dx$ for every $u\in H$ and $v\in V$,
where $\<\cpto,\cpto>$ denotes the duality pairing between $\Vp$ and~$V$.
Furthermore, we list our assumptions and notations
regarding the convex set:
\Bsist
  && \hbox{$K$ is a \revis{nonempty, closed, and} convex subset of $H$}
  \label{hpK}
  \\
  && \hbox{$\dK:H\to H$\enskip is the distance function associated to $K$}
  \label{defdK}
  \\
  && \hbox{$\PK:H\to H$\enskip is the projection operator on $K$}
  \label{defPK}
  \\
  && \hbox{$\QK:=\IH-\PK$\enskip, where $\IH:H\to H$ is the identity map}.
  \label{defQK}
\Esist
So, \accorpa{hpK}{defQK} are related to each other as follows:
for $v\in H$, we have
\Beq
  \PK v \in K 
  \aand
  \dK(v) = \normaH{\QK v} \leq \normaH{v-z}
  \quad \hbox{for every $z\in K$}.
  \label{relaz}
\Eeq
Concerning the data, \revis{we in principle} assume that
\Beq
  f \in \L2H
  \aand
  \thetaz \in V .
  \label{hpdati}
\Eeq
\Accorpa\HP hpkappa hpdati
At this point, we can state the problem \juerg{under investigation}:
given a real number $\rho>0$, we look for a triplet $(\theta,u,\sigma)$
satisfying the regularity properties
\Bsist
  && \theta \in \H1\Vp \cap \L\infty H
  \label{regtheta}
  \\
  && u \in \H1H \cap \L\infty V
  \label{regu}
  \\
  && \sigma \in \L\infty H
  \label{regsigma}
\Esist
\Accorpa\Regsoluz regtheta regsigma
and solving the problem
\Bsist
  && \< \dt\theta(t) , v >
  + \iO \nabla u(t) \cdot \nabla v
  + \rho \iO \sigma(t) \, v
  = \iO f(t) \, v
  \non
  \\
  && \quad \hbox{\aat\ and every $v\in V$}
  \label{prima}
  \\
  && u = G(\theta)
  \quad \aeQ
  \aand
  \sigma(t) \in \partial\dK(\theta(t))
  \quad \aat
  \qquad
  \vphantom\int
  \label{seconda}
  \\
  && \theta(0) = \thetaz \,.
  \label{cauchy}
\Esist
\Accorpa\Pbl prima cauchy
We observe that the system \accorpa{prima}{seconda} 
is the variational formulation of the differential inclusion \eqref{Imodif}
complemented with the no-flux \juerg{boundary} condition for $\nabla u$
(i.e.,~for $\katheta\nabla\theta$ whenever the chain rule can be applied).
Moreover, we notice that \eqref{regtheta} implies that
$\theta\in\C0\Vp$
(and even that $\theta$ is a weakly continuous $H$-valued function)
so~that the initial condition \eqref{cauchy} makes sense.
Here is our first result.

\Bthm
\label{Wellposedness}
Under the assumptions and notations \HP, suppose that $\kmin>0$.
Then, for every $\rho>0$, problem \Pbl\ has at least \juerg{one}
 solution $(\theta,u,\sigma)$
satisfying \Regsoluz\ and 
\Beq
  \theta \in \H1H \cap \L\infty V .
  \label{thetapiureg}
\Eeq
Moreover, \juerg{there is only one such solution} if $\kappa$ is a positive constant.
\Ethm

We can deal with the degenerate case $\kmin=0$
only for convex sets of obstacle or double-obstacle type.
Namely, we suppose that
\Bsist
  && \hbox{$I$\enskip is a \juerg{closed nonempty} interval}
  \label{hpI}
  \\
  && K := \graffe{v\in H:\ v(x)\in I\ \aaO} .
  \label{defK}
\Esist
\Accorpa\HPbis hpI defK
In this case, the projection on $K$ is a \juerg{pointwise} projection, i.e., for $v\in H$
\juerg{and almost every $\,x\in \Omega$, it holds that}
\Beq
  \juerg{  \hbox{$(\PK v)(x)$ is the projection of $v(x)$ on the interval $I$}}.
  \qquad
  \label{ptproj}
\Eeq
Moreover, we have to reinforce our assumptions \juerg{by postulating that}
\Beq
  f \in \L1{\Lx p}
  \aand
  \thetaz \in \Lx p
  \quad \hbox{for some $p>2$}.
  \label{hpdatibis}
\Eeq

\Bthm
\label{Wellposednessbis}
Besides the assumptions and notations \HP, suppose that \HPbis\ and \eqref{hpdatibis} hold.
Then, for every $\rho>0$, problem \Pbl\ has at least \juerg{one} solution $(\theta,u,\sigma)$
satisfying \Regsoluz\ as well~as
\Beq
  \theta \in \L\infty{\Lx p} .
  \label{regthetabis}
\Eeq
\Ethm

The main result of our paper is the next one.
It holds \revis{true} for the particular class \HPbis\ of convex subsets.
However, the degenerate case $\kmin=0$ is allowed as well.
We ensure the existence of a solution $(\theta,u,\sigma)$
whose component $\theta$ approaches 
and eventually reaches the convex set $K$ in a finite time~$\Tstar<T$,
provided that the parameter $\rho$ is large enough.
Indeed, from the statement it follows that the condition
\Beq
  \rho > \rhostar + \frac {\dK(\thetaz)} T 
  \quad \hbox{where} \quad 
  \rhostar := \norma f_{\L\infty H}
  \label{suffcond}
\Eeq
implies $\Tstar<T$.
Moreover, the speed \revis{of} approach is $\rho-\rhostar$, at least.
The precise meaning of the theorem relies on the following observation,
which follows from the regularity of $\theta$ specified by~\eqref{regtheta}:
\Beq
  \hbox{$\theta$ is an $H$-valued weakly continuous function}.
  \label{weakcont}
\Eeq
Namely, the continuous representative $\theta\in\C0\Vp$
satisfies $\theta(t)\in H$ for every $t\in[0,T]$\revis{,}
and $\theta$ is continuous from $[0,T]$ to~$H$ endowed with its weak topology.

\Bthm
\label{Goal}
Under the assumptions and notations \HP\ and \HPbis,
suppose either $\kmin>0$ or~\eqref{hpdatibis}. 
Furthermore, assume that
\Beq
  \rhostar := \norma f_{\L\infty H} < + \infty \,.
  \label{hpdatiter}
\Eeq
Then, for every $\rho>\rhostar$,
there exists a solution $(\theta,u,\sigma)$ to problem \Pbl\ with the following properties:
\Bsist
  & i)
  & \hbox{if $\thetaz\in K$,\quad then} \quad
  \hbox{$\theta(t)\in K$ for every $t\in[0,T]$} 
  \label{trivial}
  \\
  & ii)
  & \hbox{if $\thetaz\not\in K$,\quad there exists $\Tstar\in(0,T]$ satisfying} \quad
  \Tstar \leq \frac{\dK(\thetaz)}{\rho-\rhostar}
  \quad \hbox{such that}
  \non
  \\
  \noalign{\vskip -3pt}
  && \frac d{dt} \, \dK(\theta(t)) \leq - (\rho-\rhostar)
  \quad \hbox{in the sense of distributions on $(0,\Tstar)$}
  \label{decreasing}
  \\
  \noalign{\vskip 4pt}
  && \theta(t) \in K
  \quad \hbox{for every $t\in[0,T]$ such that} \quad t > \Tstar .
  \label{goal}
\Esist
In particular, in the case $ii)$, the function
$t\mapsto\dK(\theta(t))$ is strictly decreasing on~$[0,\Tstar]$.

\Ethm

We close this section with a list of \revis{denotations} and tools.
Throughout the paper, 
$\norma\cpto_X$ denotes the norm in the generic Banach space~$X$
or in a power \revis{thereof}.
However, we simply write $\norma\cpto_p$ for the standard norm in~$\Lx p$.
Moreover, we repeatedly use the \juerg{denotation}
\Beq
  Q_t := \Omega \times (0,t)
  \quad \hbox{for $0<t\leq T$}
  \label{defQt}
\Eeq
as well as the Young inequalities
\Bsist
  && ab\leq \delta a^2 + \frac 1 {4\delta}\,b^2
  \aand
  ab \leq \badtheta a^{\frac 1\badtheta} + (1-\badtheta) b^{\frac 1{1-\badtheta}}
  \non
  \\
  && \quad \hbox{for every $a,b\geq 0$, \ $\delta>0$, \ and \ $\badtheta\in(0,1)$}
  \label{young}
\Esist
and \juerg{H\"older's} inequality.
Furthermore, we account for the compact embedding $V\subset H$.
Finally, we follow a general rule to denote constants\revis{: the}
small-case symbol $c$ stands for different constants which depend only
on~$\Omega$, on the final time~$T$, the structure of the problem
and on the constants and the norms of
the functions involved in the assumptions of our statements.
A symbol like $c_\delta$ signals that the constant can depend also on the parameter~$\delta$.
Hence, the meaning of~$c$ (or~$c_\delta$) might
change from line to line and even \revis{within} the same chain of equalities or inequalities.


\section{Approximation}
\label{APPROX}
\setcounter{equation}{0}

In this section, we introduce an approximating problem 
which depends on the parameters $\eps,\alpha\in(0,1)$
and is useful to establish some parts of our results.
We \juerg{could have} decided to take $\alpha=\eps$ to reach the same goal.
However, we think that the choice of two different parameters
could \juerg{prove to} be more suitable for the numerical \juerg{treatment}.
From one side, we replace the function $\kappa$ 
by a strictly positive $\keps$ in order to ensure uniform parabolicity.
On the other hand, we \regulariz e the subdifferential~$\partial\dK$.
However, for the sake of simplicity,
we often avoid stressing the dependence on both parameters in the notation
and write, e.g.,~$\thetaeps$ instead of~$\theta_{\eps,\alpha}$.
We introduce the functions $\keps$ and $\Geps$ as follows\revis{:}
\Bsist
  && \keps := \kappa
  \quad \hbox{if $\kmin>0$\revis{,}}
  \aand
  \keps := \kappa + \alpha
  \quad \hbox{if $\kmin=0$\revis{,}}
  \label{defkeps}
  \\
  && \Geps(r) := \int_0^r \keps(s) \, ds
  \quad \hbox{for $r\in\erre$}.
  \label{defGeps}
\Esist
Moreover, let $\dKeps:H\to\erre$ and $D\dKeps:H\to H$
be the Moreau--Yosida \regulariz ations
of the nondifferentiable function $\dK$ and of its subdifferential~$\partial\dK$.
Thus, for $v\in H$, we~have \juerg{that}
\Bsist
  && \dKeps(v) := \inf_{z\in H} \Bigl( \dK(z) + \sfrac 1{2\eps} \, \normaH{z-v}^2 \Bigr)
  \label{genmoreau}
  \\
  && \hbox{$D\dKeps(v)$ is the gradient of $\dKeps$ at $v$, i.e., the unique element of $\partial\dKeps(v)$}.
  \qquad
  \label{genyosida}
\Esist
The \juerg{statement} \eqref{genyosida} means that the map
$H\ni z\mapsto(D\dKeps(v),z)_H$ is the \frechet\ derivative of $\dKeps$ at~$v$
($\dKeps$~is \frechet\ differentiable, indeed).
We recall that the subdifferential~$\partial\dKeps$, 
which we identify with the single-valued map~$D\dKeps$,
actually is the Yosida \regulariz ation of~$\partial\dK$, thus \Lip\ continuous
with \Lip\ constant~$1/\eps$
(see, e.g., \cite[p.~28 and Prop.~2.11, p.~39]{Brezis}).
\juerg{These maps can be given} explicitly, as shown in the next lemma.
As we could not find precise references on it,
we proved the result in the Appendix.

\Blem
\label{MorYos}
Let $H$ be a Hilbert space.
With the assumptions and notations \accorpa{defdK}{defQK} and \accorpa{genmoreau}{genyosida}, the formulas\pier{%
\Beq
  D\dKeps(v) = \frac {\QK v} {\max\{\eps,\dK(v)\}}
  \label{moryos1}
\Eeq
and 
\Beq  
\dKeps(v) = \int_0^{\dK(v)} \min\{s/\eps,1\} \, ds
  \label{moryos2}
\Eeq
}%
hold true for every $v\in H$.
\Elem

At this point, we introduce the approximating problem.
It consists in finding a triplet $(\thetaeps,\ueps,\sigmaeps)$ satisfying
\Beq
  \thetaeps, \ueps \in \H1H \cap \L\infty V
  \aand
  \sigmaeps \in \L\infty H
  \label{regsoluzeps}
\Eeq
and solving the variational problem
\Bsist
  && \iO \dt\thetaeps(t) \, v 
  + \iO \nabla\ueps(t) \cdot \nabla v
  + \rho \iO \sigmaeps(t) \, v
  = \iO f(t) \, v
  \non
  \\
  && \quad \hbox{\aat\ and every $v\in V$}
  \label{primaeps}
  \\
  && \ueps = \Geps(\thetaeps)
  \quad \aeQ
  \aand
  \sigmaeps(t) = D\dKeps(\thetaeps(t))
  \quad \aat
  \qquad
  \vphantom\int
  \label{secondaeps}
  \\
  && \thetaeps(0) = \thetaz \,.
  \label{cauchyeps}
\Esist
\Accorpa\Pbleps primaeps cauchyeps

\Bthm
\label{Welleps}
Under the assumptions and notations \HP\ and \accorpa{defkeps}{genyosida},
for every $\eps,\alpha\in(0,1)$,
problem \Pbleps\ has a unique solution $(\thetaeps,\ueps,\sigmaeps)$ satisfying~\eqref{regsoluzeps}.
\Ethm

The rest of the section is devoted to prove \juerg{this} well-posedness result.
We first establish the existence of a solution \revis{via} a fixed point argument.
Concerning the symbols $\theta$, $u$ and $\sigma$ we often use,
we point out that they have nothing to do with the original problem~\Pbl,
which is out of interest at the moment.

\step
Existence

For a given $\thetabar\in\LQ2$, 
we look for a solution $\theta\in\H1H\cap\L\infty V$ 
to the problem
\Bsist
  && \iO \dt\theta(t) \, v 
  + \iO \keps(\theta) \nabla\theta(t) \cdot \nabla v
  = \iO \bigl( f(t) - \rho \, \sigma(t) \bigr) \, v
  \non
  \\
  && \quad \hbox{\aat\ and every $v\in V$}
  \label{primabar}
  \\
  \noalign{\smallskip}
  && \theta(0) = \thetaz \,.
  \label{cauchybar}
\Esist
\Accorpa\Pblbar primabar cauchybar
where $\sigma(t) := D\dKeps(\thetabar(t))$ \aat .
As $\keps$ is a continuous function such that $\alpha\leq\keps\leq\kmax+1$, 
problem \Pblbar\ has a unique solution $\theta$ satisfying the prescribed regularity.
Moreover, by testing \eqref{primabar} with $v=\theta(t)$\revis{,}
and noting that $\normaH{\sigma(t)}\leq1$ by \pier{\eqref{moryos1}},
we immediately obtain \revis{that}
\Bsist
  && \frac 12 \, \normaH{\theta(t)}^2
  \leq \juerg{\frac 12 \|\vartheta_0\|_H^2 \,+}\iot 
	\bigl( \normaH{f(s)}+\rho \bigr) \normaH{\theta(s)} \, ds 
  \non
  \\
  && \leq \frac 12 \, \normaH\thetaz^2
  + \norma f_{\L2H}^2 + \rho^2 T
  + \frac 12 \iot \normaH{\theta(s)}^2 \, ds .
  \non
\Esist
By applying the Gronwall lemma, we deduce a bound in $\L\infty H$ and infer that
\Beq
  \norma\theta_{\L2H} \leq R
  \label{defR}
\Eeq
for some constant $R$ depending only on the data, $T$ and~$\rho$.
At this point, we denote by $\calB^2$ and $\calB^\infty$ the closed unit balls 
of $\L2H$ and $\L\infty H$, respectively,
set $\calK:=R\,\calB^2$ and define the maps
$\calS:\calK\to\calB^\infty$ and $\calF:\calK\to\calK$ by setting for $\thetabar\in\calK$
\sl 
\begin{align}
  & \calS(\thetabar) := \sigma
  \quad \hbox{given by} \quad
  \sigma(t) := D\dKeps(\thetabar(t))
  \quad \aat 
  \label{defsigma}
  \\
  & \juerg{\hbox{$\theta:=\calF(\thetabar)\in\H1H\cap\L\infty V$ is the 
	unique solution to \Pblbar}}.
  \label{defF}
\end{align}
\rm
We verify that we can apply the Schauder fixed point theorem to~$\calF$
with respect to the strong topology of~$\L2H$.
Clearly, $\calK$~is \juerg{nonempty, bounded, convex and closed}.
Next, if $\thetabar\in\calK$ and $\theta:=\calF(\thetabar)$,
then $\theta$, $u:=\Geps(\theta)$ and $\sigma:=\calS(\thetabar)$ satisfy
$\dt\theta - \Delta u = f - \rho \, \sigma$
in the sense of distributions on~$Q$, in principle, 
then \aeQ\ since $\dt\theta$ and the \rhs\ belong to~$\LQ2$.
Moreover, $u$ satisfies the homogeneous Neumann boundary condition.
By multiplying the above equation by $\dt u=\dt\theta/\keps(\theta)$
and also recalling that $\keps\leq\kmax+1$ and that $\norma\sigma_{\L\infty H}\leq1$,
we easily obtain
\Beq
  \norma{\dt\theta}_{\L2H} + \norma{\nabla u}_{\L\infty H} \leq c \,,
  \quad \hbox{whence also} \quad
  \norma{\nabla\theta}_{\L\infty H} \leq \ca
  \non
\Eeq
since $\nabla\theta=\nabla u/\keps(\theta)$ and $\keps\geq\alpha$.
We conclude that
\Beq
  \norma\theta_{\H1H\cap\L\infty V} \leq \ca \,.
  \label{stimaF}
\Eeq
By the Aubin-Lions lemma (see, e.g., \cite[Thm.~5.1, p.~58]{Lions}), 
we see that $\calF(\calK)$ is relatively compact \juerg{in $L^2(0,T;H)$}.
Finally, we check that $\calF$ is continuous.
Let $\thetabarn,\thetabar\in\calK$ be such that $\thetabarn\to\thetabar$ strongly in~$\LQ2$,
and set $\sigman:=\calS(\thetabarn)$ and $\thetan:=\calF(\thetabarn)$.
Then, $\sigman$ converge to $\sigma:=\calS(\thetabar)$ strongly in $\L2H$ 
since $D\dKeps$ is \Lip\ continuous on~$H$.
Furthermore, estimate \eqref{stimaF} holds for~$\thetan$.
Therefore, we have (for a subsequence)
\Bsist
  & \thetan \to \theta
  & \quad \hbox{weakly star in $\H1H\cap\L\infty V$}
  \non
  \\
  && \quad \hbox{strongly in $\L2H$ and \aeQ}
  \non
\Esist
for some $\theta\in\H1H\cap\L\infty V$ which necessarily belongs to~$\calK$.
Since $\keps$ is continuous and bounded,
we infer that $\keps(\thetan)$ converges to $\keps(\theta)$ strongly in $\LQ p$
for every $p\in[1,+\infty)$.
Thus, it is \sfw\ to deduce that $\theta$ solves~\Pblbar, i.e., that $\theta=\calF(\thetabar)$,
and that the convergence holds for the whole sequence~$\graffe{\thetan}$.
Therefore, $\calF$~is continuous and we conclude that it has at least a fixed point.
Now, if $\thetaeps$ is a fixed point of~$\calF$ and we set
$\ueps:=\Geps(\thetaeps)$ and $\sigmaeps:=D\dKeps(\thetaeps)$,
one easily sees that the triplet $(\thetaeps,\ueps,\sigmaeps)$ satisfies~\eqref{regsoluzeps}
and it is clear that it is a solution to problem \Pbleps.

\step
Uniqueness

Let $(\theta_1,u_1,\sigma_1)$ and $(\theta_2,u_2,\sigma_2)$ be two solutions of problem \Pbleps\ 
satisfying the regularity requirement \eqref{regsoluzeps}.
We write \eqref{primaeps} for both of them, take the difference, and integrate with respect to time.
We have for almost every $s\in(0,T)$ and every $v\in V$
\Beq
  \iO \bigl( \theta_1(s) - \theta_2(s) \bigr) \, v
  + \iO \nabla \bigl( 1*(u_1-u_2) \bigr)(s) \cdot \nabla v
  = - \rho \iO \bigl( 1*(\sigma_1-\sigma_2) \bigr)(s) \, v
  \non
\Eeq
with the general notation $(1*v)(s):=\int_0^s v(\tau)\,d\tau$.
Now, we choose $v=(u_1-u_2)(s)$ and integrate over $(0,t)$ with respect to~$s$.
We obtain
\Beq
  \intQt (\theta_1-\theta_2) (u_1-u_2)
  + \frac 12 \iO |\bigl( 1*\nabla(u_1-u_2) \bigr)(t)|^2
  = - \rho \intQt \bigl( 1*(\sigma_1-\sigma_2) \bigr) (u_1-u_2) .
  \non
\Eeq
We recall that $\alpha\leq\Geps'\leq\kmax+1$
and ignore the second term on the \lhs, which is nonnegative. 
Furthermore, we owe to the Young inequality.
We deduce that
\Beq
  \alpha \iot \normaH{(\theta_1-\theta_2)(s)}^2 \, ds
  \leq \frac \alpha 2 \iot \normaH{(\theta_1-\theta_2)(s)}^2 \, ds
  + \frac c\alpha \intQt |1*(\sigma_1-\sigma_2)|^2 .
  \label{peruniceps}
\Eeq
Now, we use the \Holder\ inequality
and account for the $(1/\eps)$-\Lip\ continuity of~$D\dKeps$
(as~a map from $H$ into itself).
We have for every $s\in[0,T]$
\Bsist
  && \normaH{\bigl( 1*(\sigma_1-\sigma_2) \bigr)(s)}^2
  = \normaH{{\textstyle\int_0^s(\sigma_1-\sigma_2)(\tau)\,d\tau}}^2
  \non
  \\
  && \leq c \int_0^s \normaH{(\sigma_1-\sigma_2)(\tau)}^2 \, d\tau
  \leq \ceps \int_0^s \normaH{(\theta_1-\theta_2)(\tau)}^2 \, d\tau
  \non
\Esist
and deduce that
\Beq
  \intQt |1*(\sigma_1-\sigma_2)|^2
  \leq \ceps \iot \Bigl( \int_0^s \normaH{(\theta_1-\theta_2)(\tau)}^2 \, d\tau \Bigr) \, ds .
  \non
\Eeq
Coming back to \eqref{peruniceps} and applying the Gronwall lemma,
we conclude that $\theta_1=\theta_2$, whence also $u_1=u_2$ and $\sigma_1=\sigma_2$.


\section{Well-posedness}
\label{PROOFS}
\setcounter{equation}{0}

This section deals with Theorems~\ref{Wellposedness} and~\ref{Wellposednessbis}.
In order to prove the statements regarding existence,
we start from the solution $(\thetaeps,\ueps,\sigmaeps)$ to the approximating problem \Pbleps\
and perform a number of a~priori estimates \juerg{in which all of the occurring constants
$c>0$ will be independent of both $\varepsilon$ and $\alpha$}.

\step
First a priori estimate

We test \eqref{primaeps} by $\thetaeps$ and have
\Beq
  \frac 12 \iO |\thetaeps(t)|^2
  + \intQt \nabla\ueps \cdot \nabla\thetaeps
  = \frac 12 \iO |\thetaz|^2
  + \intQt \bigl( f - \rho\,\sigmaeps \bigr) \thetaeps \,.
  \non
\Eeq
We recall that $\norma{\sigmaeps}_{\L\infty H}\leq1$ (cf.~\pier{\eqref{moryos1}})
and observe that there hold, \aeQ,
\Bsist
  && \pier{|\ueps| 
  = |\Geps(\thetaeps)|
  \leq (\kmax+1) |\thetaeps|}
  \non
  \\
  && \pier{\nabla\ueps \cdot \nabla\thetaeps
  = \frac 1 {\keps(\thetaeps)} \, |\nabla\ueps|^2
  \geq \frac 1 {\kmax + 1} \, |\nabla\ueps|^2}
  \non
  \\
  \noalign{\smallskip}
  && \nabla\ueps \cdot \nabla\thetaeps
  = \keps(\thetaeps) |\nabla\thetaeps|^2
  \geq \kmin |\nabla\thetaeps|^2 .
  \non
\Esist
Hence, by also owing to the Gronwall lemma, we easily deduce that
\Bsist
  && \norma\thetaeps_{\L\infty H}
  + \norma\ueps_{\L\infty H\cap\L2V}
  \leq c 
  \label{primastima}
  \\
  && \norma\thetaeps_{\L2V} \leq c
  \quad \hbox{provided that $\kmin>0$}.
  \label{primastimabis}
\Esist
By comparison in the variational equation \eqref{prima}, \juerg{we infer from 
\eqref{primastima}} that
\Beq
  \norma{\dt\thetaeps}_{\L2\Vp} \leq c \,.
  \label{stimadt}
\Eeq

\step 
Second a priori estimate

We notice that \eqref{prima} can be written as
\Beq
  \dt\thetaeps - \div \bigl( \keps(\thetaeps) \nabla\thetaeps \bigr)
  = f - \rho \, \sigmaeps
  \quad \aeQ \revis{,}
  \non
\Eeq
with the homogeneous Neumann boundary condition for~$\ueps$.
By multiplying by $\dt\ueps$, integrating over~$Q_t$, and applying the \Holder\ and Young inequalities, we obtain
\Beq
  \intQt \frac 1 {\keps(\thetaeps)} \, |\dt\ueps|^2
  + \frac 12 \iO |\nabla\ueps(t)|^2
  \leq c + \frac 1 {2(\kmax+1)} \intQt |\dt\ueps|^2 .
  \non
\Eeq
As $\keps(\thetaeps)\leq\kmax+1$, we deduce on account of~\eqref{primastima}
\revis{that}
\Beq
  \norma\ueps_{\H1H\cap\L\infty V} \leq c \,.
  \label{secondastima}
\Eeq
Moreover, from the identities
\Beq
  \dt\thetaeps = \frac {\dt\ueps} {\keps(\thetaeps)}
  \aand
  \nabla\thetaeps = \frac {\nabla\ueps} {\keps(\thetaeps)}
  \non
\Eeq
we infer that
\Beq
  \norma\thetaeps_{\H1H\cap\L\infty V} \leq c
  \quad \hbox{provided that $\kmin>0$}.
  \label{stimadtbis}
\Eeq

\step
Convergence

At this point, by standard compactness results
(in particular, we owe to the Aubin-Lions lemma proved, e.g., 
in \cite[Thm.~5.1, p.~58]{Lions})), 
we deduce that
\Bsist
  & \thetaeps \to \theta
  & \quad \hbox{weakly star in $\H1\Vp\cap\L\infty H$}
  \non
  \\
  & 
  & \quad \hbox{and strongly in $\L2\Vp$}  
  \label{convtheta}
  \\
  & \ueps \to u
  & \quad \hbox{weakly star in $\H1H\cap\L\infty V$}
  \non
  \\
  & 
  & \quad \hbox{strongly in $\L2H$ and \aeQ}  
  \label{convu}
  \\
  & \sigmaeps \to \sigma
  & \quad \hbox{weakly star in $\L\infty H$}
  \label{convsigma}
\Esist
for some triplet $(\theta,u,\sigma)$, 
as $(\eps,\alpha)$ tends to~$(0,0)$,
at least for a subsequence.
Moreover, we also have
\Bsist
  & \thetaeps \to \theta
  & \quad \hbox{weakly star in $\H1H\cap\L\infty V$}
  \non
  \\
  & 
  & \quad \hbox{strongly in $\L2H$ and \aeQ}  
  \quad \hbox{provided that $\kmin>0$}.
  \label{convthetabis}
\Esist
Thus, it is clear that $(\theta,u,\sigma)$ is a solution to \Pbl\ 
satisfying the regularity requirements 
stated in Theorems~\ref{Wellposedness} and~\ref{Wellposednessbis}\juerg{, with the 
exception of \eqref{regthetabis}},
whenever we prove that
$u=G(\theta)$ \aeQ\ and that $\sigma(t)\in\partial\dK(\theta(t))$ \aat.
At this point, we have to distinguish the different cases 
corresponding to the above statements.

\step
Conclusion of the existence proof in the uniformly parabolic case

We complete the proof of the existence part of Theorem~\ref{Wellposedness} by assuming $\kmin>0$.
Thus, $\keps=\kappa$ and $\Geps=G$.
Thanks to the pointwise convergence~(a.e.)\ given by \eqref{convu} and~\eqref{convthetabis}
and to the continuity of~$G$,
we immediately deduce that $u=G(\theta)$ \aeQ.
As for the second \pier{condition in}~\eqref{seconda}, we owe to the strong convergence~\eqref{convthetabis}
and apply, e.g., \cite[Lemma~2.3, p.~38]{Barbu} 
to the maximal monotone operator induced on $\L2H$ by~$\partial\dK$.

On the contrary, for the degenerate case allowed in Theorem~\ref{Wellposednessbis}, 
some more work has to be done.
Concerning the relation $u=G(\theta)$ that we have to prove, we observe that
$\Geps(r)\geq G(r)$ for $r\geq0$ and $\Geps(r)\leq G(r)$ for $r\leq0$ since $\kappa\leq\keps$.
It follows that $\inf\Geps\pier{{}= -\infty } \leq\inf G$ and $\sup\Geps\pier{{}= +\infty }\geq\sup G$, i.e.,
the domain $D(\Geps^{-1})$ of $\Geps^{-1}$ includes the domain $D(G^{-1})$ of~$G^{-1}$.

\Blem
\label{ConvinvG}
The following convergence holds true\revis{:}
\Beq
  \Geps^{-1}(s) \to G^{-1}(s) \quad
  \hbox{uniformly in every compact subset of $D(G^{-1})$}.
  \label{convinvG}
\Eeq
\Elem

\Bdim
We first establish the pointwise convergence.
This trivially holds if $s=0$.
Assume $s>0$.
Then, $\Geps^{-1}(s)>G_{\alpha'}^{-1}(s)>G^{-1}(s)$ for $\alpha'\in(0,\alpha)$
since $G(r)<G_{\alpha'}(r)<\Geps(r)$ for every $r>0$.
Thus the limit $\ell$ of $\Geps^{-1}(s)$ as $\alpha\searrow0$ exists and satisfies
$\ell\geq\ell_0:=G^{-1}(s)$.
It follows that the constant $s=\Geps(\Geps^{-1}(s))$ converges to~$G(\ell)$,
i.e., that $s=G(\ell)$.
As $G(\ell_0)=s$ and $G$ is one-to-one by~\eqref{hpG}, we conclude that $\ell=\ell_0$.
Thus, we have proved that $\Geps^{-1}(s)$ converges to $G^{-1}(s)$ pointwise
in the interval $D(G^{-1})\cap[0,+\infty)$.
Since the convergence is monotone and the limit $G^{-1}$ is continuous,
the convergence is uniform on every compact subset (Dini's theorem).
As the case of negative values of $s$ is similar, \eqref{convinvG} is proved.
\Edim

At this point\revis{,} we can go on and show that $u=G(\theta)$ \aeQ.
From Lemma~\ref{ConvinvG} and the pointwise convergence \eqref{convu} of $\ueps$ to $u$
we infer that $\thetaeps=\Geps^{-1}(\ueps)$ converges to $G^{-1}(u)$ \aeQ.
Since $\thetaeps$ converges to $\theta$ weakly in $\LQ2$ by~\eqref{convtheta},
we conclude that \pier{(see, e.g., \cite[Lemme~1.3, p.~12]{Lions})}
$G^{-1}(u)=\theta$, i.e., $u=G(\theta)$, \aeQ.
As a by-product, there holds the convergence
\Beq
  \thetaeps \to \theta
  \quad \aeQ \revis{,}
  \label{convae}
\Eeq
and we use it to prove that $\sigma(t)\in\partial\dK(\theta(t))$ \aat.
Here\revis{,} we owe to assumptions \HPbis\ and \eqref{hpdatibis}
on the convex~$K$ and on the data.
For convenience\revis{,} we set
\Beq
  \PI : \erre \to \erre \quad
  \hbox{is the projection on $I$}
  \aand
  \QI := I_{\erre} - \PI \revis{,}
  \non
\Eeq
where $I_{\erre}:\erre\to\erre$ is the identity map,
so that \eqref{ptproj} reads
$(\PK v)(x)=\PI(v(x))$ \aeO\revis{,}\ for every $v\in H$.
Hence, \pier{\eqref{moryos1}} becomes
\Beq
  \bigl( D\dKeps(v) \bigr)(x)
  = \frac {\QI \pier{(} v(x) \pier{)}} {\max\{\eps,\dK(v)\}}
  \quad \aeO, \quad
  \hbox{for every $v\in H$}.
  \label{ptyos}
\Eeq

\step
A new a priori estimate

We define the truncation operator $\Tn:\erre\to\erre$ by setting
$\Tn(r):=\max\{-n,\min\{r,n\}\}$ for $r\in\erre$
and use the notation $(r)^q:=|r|^q\sign r$ (with $\sign0=0$) for $r\in\erre$ and $q>0$.
Next, we take $\thetastar\in I$,
set $\zeta:=\thetaeps-\thetastar$ for convenience
and test \eqref{prima} by $(\Tn(\zeta))^{p-1}$.
By integrating over~$Q_t$ and recalling that $p>2$, we obtain that
\Bsist
  && \iO \Tnp \bigl( \zeta(t) \bigr)
  + (p-1) \intQt \keps(\thetaeps) \, |\Tn(\zeta)|^{p-2} \, \Tn'(\zeta) \, |\nabla\thetaeps|^2
  + \rho \intQt \sigmaeps \, (\Tn(\zeta))^{p-1}
  \qquad
  \non
  \\
  && = \iO \Tnp \bigl( \zeta(0) \bigr)
  + \intQt f \, (\Tn(\zeta))^{p-1}\,,
  \label{pernew}
\Esist
where we have set
\Beq
  \Tnp(r) := \int_0^r (\Tn(s))^{p-1} \, ds
  \quad \hbox{for $r\in\erre$}.
  \non
\Eeq
We observe that $p\,\Tnp(r)\geq|\Tn(r)|^p$ for every $r\in\erre$.
Indeed, if $0\leq r\leq n$, we have $p\,\Tnp(r)=r^p=|\Tn(r)|^p$;
if $r>n$, then $p\,\Tnp(r)\geq p\,\Tnp(n)=n^p=|\Tn(r)|^p$.
On the other hand, both $\Tnp$ and $|\Tn|^p$ are even functions.
Therefore, we have
\Beq
  \iO \Tnp \bigl( \zeta(t) \bigr)
  \geq \frac 1p \, \norma{\Tn(\zeta(t))}_p^p \,.
  \non
\Eeq
The second term of \eqref{pernew} is nonnegative.
For the third one, we note that the following pairs of functions share their \juerg{signs}: 
$\QI\pier{(}\thetaeps \pier{)}$ and $\zeta$ since $\thetastar\in I$;
$\sigmaeps$ and $\QI\thetaeps$ thanks to \eqref{ptyos};
$(\Tn(\zeta))^{p-1}$ and $\zeta$ by our definition of~$(\cpto)^{p-1}$.
Thus, the same holds for $\sigmaeps$ and $(\Tn(\zeta))^{p-1}$\revis{,}
so that their product is nonnegative.
To deal with the \rhs, we use assumption \eqref{hpdatibis} on $\thetaz$ and~$f$.
First, we notice that $\Tnp(r)\leq|r|^p$ for every $r\in\erre$
so that the first integral can be estimated from above by $\norma{\thetaz-\thetastar}_p^p$.
Next, we account for the \Holder\ inequality and treat the last term as follows:
\Bsist
  && \intQt f \, (\Tn(\zeta))^{p-1}
  \leq \iot \norma{f(s)}_p \, \norma{(\Tn(\zeta(s)))^{p-1}}_{p'}
  \non
  \\
  && = \iot \norma{f(s)}_p \, \bigl( \norma{\Tn(\zeta(s))}_p^p \bigr)^{1/p'} \, ds
  \leq c \iot \norma{f(s)}_p \, \bigl( 1 + \norma{\Tn(\zeta(s))}_p^p \bigr) \, ds .
  \non
\Esist
Therefore, coming back to \eqref{pernew}
and recalling that $f\in\L1{\Lx p}$,
we can apply the Gronwall lemma and deduce that $\norma{\Tn(\zeta)}_{\L\infty{\Lx p}}\leq c$,
whence immediately
\Beq
  \norma\thetaeps_{\L\infty{\Lx p}} \leq c \,.
  \label{stimaLp}
\Eeq

\step
Conclusion of the proof of Theorem~\ref{Wellposednessbis}

Estimate \eqref{stimaLp} ensures the further regularity \eqref{regthetabis}
for the function $\theta$ given by~\eqref{convtheta}.
On the other hand, the pointwise convergence~\eqref{convae}, our assumption $p>2$\pier{, \eqref{stimaLp}\revis{,} and the Egorov theorem imply that}
\Beq
  \thetaeps \to \theta
  \quad \hbox{strongly in~$\LQ2$}.
  \label{strongtheta}
\Eeq
Thus, the strong convergence \eqref{convthetabis} 
already established in the uniformly parabolic case holds also in the present one.
Therefore, we can combine it with the weak convergence of $\sigmaeps$ to $\sigma$ in $\LQ2$
ensured by \eqref{convsigma} and proceed as before, i.e.,
we apply, e.g., \cite[Lemma~2.3, p.~38]{Barbu} 
to the maximal monotone operator induced on $\L2H$ by~$\partial\dK$.
We conclude that $\sigma(t)\in\partial\dK(\theta(t))$ \aat,
and the proof of Theorem~\ref{Wellposednessbis} is complete.\revis{\qed}

\step
Uniqueness

We prove the last sentence of Theorem~\ref{Wellposedness}
by assuming that $\kappa$ is a positive constant.
We pick two solutions $(\theta_i,u_i,\sigma_i)$, $i=1,2$,
write \eqref{prima} for both of them,
and test the difference by $\theta_1-\theta_2$.
We obtain
\Beq
  \frac 12 \iO |(\theta_1-\theta_2)(t)|^2
  + \kappa \intQt |\nabla(\theta_1-\theta_2)|^2
  + \rho \intQt (\sigma_1-\sigma_2)(\theta_1-\theta_2)
  = 0 .
  \non
\Eeq
\juerg{All of} the terms on the \lhs\ are nonnegative,
the third one since $\sigma_i(t)\in\partial\dK(\theta_i(t))$ \aat\ 
and $\partial\dK:H\to2^H$ is monotone.
We immediately deduce that $\theta_1=\theta_2$, whence also $u_1=u_2$.
By comparison in~\eqref{prima}, we infer that $\sigma_1=\sigma_2$ as 
well.\revis{\qed}

\section{Proof of Theorem~\ref{Goal}}
\label{GOAL}
\setcounter{equation}{0}

It suffices to prove that the solution $(\theta,u,\sigma)$ constructed in the previous section
satisfies the conditions of the statement.
Therefore, we keep the notation already used
and the boundedness and convergence specified in the above proofs
(the latter for a not relabeled subsequence $(\eps,\alpha)\to(0,0)$, as usual).
In~fact, we only need to know that $(\theta,u,\sigma)$
is a solution to \Pbl\ satisfying the regularity conditions \Regsoluz\
and to account~for
\Beq
  \norma\thetaeps_{\L\infty H} \leq c
  \aand
  \thetaeps \to \theta
  \quad \hbox{strongly in $\LQ2$}.
  \label{strong}
\Eeq

\step
The differential inequality

We test \eqref{primaeps} by $\sigmaeps(t)=D\dKeps(\thetaeps(t))$, integrate over $\Omega$\revis{,} and obtain\revis{,}~\aat\revis{,}
\Beq
  \frac d{dt} \, \dKeps(\thetaeps(t))
  + \iO \keps(\thetaeps(t)) \nabla\thetaeps(t) \cdot \nabla\sigmaeps(t)
  + \rho \normaH{\sigmaeps(t)}^2
  = \iO f(t) \, \sigmaeps(t) \,.
  \label{perbasic}
\Eeq
We observe that the second term on the \lhs\ is nonnegative
\pier{on account of}~\eqref{ptyos}.
\revis{Indeed,} for almost all $x\in\Omega$\revis{, we have}
\Beq
  \sigmaeps(x,t) = \frac {\QI \pier{(} \thetaeps(x,t)\pier{)}} {\max\{\eps,\dK(\pier{\thetaeps}(t))\}}\,\revis{,}
  \quad \hbox{whence} \quad
  \nabla\sigmaeps(x,t) = \frac {\QI'(\thetaeps(x,t))\nabla\thetaeps(x,t)} {\max\{\eps,\dK(\pier{\thetaeps}(t))\}} \,.
  \non
\Eeq
On the other hand, both $\keps$ and $\QI'$ are nonnegative functions.
As for the \rhs, we recall assumption \eqref{hpdatiter} on~$f$
and that $\normaH{\sigmaeps(t)}\leq1$ \aat\ by \pier{\eqref{secondaeps} and~\eqref{moryos1}}.
Therefore, we deduce from~\eqref{perbasic}~\revis{that}
\Beq
  \frac d{dt} \, \dKeps(\thetaeps(t))
  + \rho \normaH{\sigmaeps(t)}^2 
  \leq \rhostar := \norma f_{\L\infty H} .
  \label{basic}
\Eeq
Notice that the definition of $\rhostar$ agrees with~\eqref{hpdatiter}.
We observe that the formulas \pier{\eqref{moryos1}--\eqref{moryos2}} imply that
$\normaH{D\dKeps(v)}=1$ if $\dK(v)>\eps$
and that $\dK(v)>\eps$ if and only if $\dKeps(v)>\eps/2$.
Hence, we deduce from \eqref{basic} that the functions 
$\psieps\in W^{1,1}(0,T)$ and $\phieps\in L^\infty(0,T)$\revis{,}
defined~by
\Beq
  \psieps(t) := \dKeps(\thetaeps(t))
  \aand
  \phieps(t) := \normaH{\sigmaeps(t)}^2 \revis{,}
  \label{defpsiphi}
\Eeq
satisfy
\Bsist
  && \psieps'(t) + \rho \, \phieps(t) \leq \rhostar
  \quad \aat
  \label{diffineq}
  \\
  && \phieps(t) = 1
  \quad \hbox{a.e.\ in the set where $\psieps>\eps/2$}.
  \label{good}
\Esist

\Blem
\label{Diffineq}
Let $\psi\in W^{1,1}(0,T)$ and $\phi\in L^\infty(0,T)$
and assume that for some positive constant $\rho,\delta,\gamma$
the following conditions hold
\Bsist
  && \psi'(t) + \rho \, \phi(t) \leq \rho-\delta
  \quad \aat
  \non
  \\
  && \phi(t) = 1
  \quad \hbox{a.e.\ in the set where $\psi>\gamma$}.
  \non
\Esist
Then, we have $\psi(t)\leq\gamma$ for every $t\in[0,T]$ if $\psi(0)\leq\gamma$.
Moreover, if $\psi(0)>\gamma$, \revis{then there exists some} 
$\Tg\in(0,T]$ satisfying
$\Tg \leq \bigl(\psi(0)-\gamma\bigr)/\delta$ such that
\Beq
  \psi'\leq - \delta
  \quad \hbox{a.e.\ in $(0,\Tg)$}
  \aand
  \psi(t) \leq \gamma
  \quad \hbox{for every $t\in[0,T]$ such that} \quad t > \Tg \,.
  \label{tesi}
\Eeq
\Elem

\Bdim
We start with the following statement:
\Beq
  \hbox{if $\tz\in[0,T)$ and $\psi(\tz)\leq\gamma$,
    then $\psi(t)\leq\gamma$ for every $t\in[\tz,T]$}.
  \label{prelim}
\Eeq
By contradiction, assume that there exists some $t\in(\tz,T]$ such that $\psi(t)>\gamma$ 
and consider the set $A:=\graffe{s\in(\tz,t]:\ \psi(\tau)>\gamma\ \hbox{for every}\ \tau\in (s,t]}$.
Then, $A$~is \revis{nonempty} and $\tu:=\inf A$ satisfies
$\tu\geq\tz$, whence $\psi(\tu)=\gamma$.
Moreover, $\tu<t$ and $\psi(\tau)>\gamma$ for every $\tau\in(\tu,t]$.
From the assumptions\revis{,} it follows that $\psi'\leq-\delta$ a.e.\ in~$(\tu,t)$.
Thus, $\psi(t)\leq\psi(\tu)-\delta(t-\tu)<\gamma$, a~contradiction.
Therefore, \eqref{prelim} is established.
From \eqref{prelim} we deduce the first sentence of the statement
regarding the case $\psi(0)\leq\gamma$.
Assume now \revis{that $\psi(0)>\gamma$,}
and set $\Tg:=\sup\graffe{t\in(0,T]:\ \psi(s)>\gamma\ \hbox{for every}\ s\in[0,t)}$.
Then, $\Tg>0$ and $\psi(s)>\gamma$ for every $s\in[0,\Tg)$.
Our assumptions imply \revis{that} $\phi=1$ and $\psi'\leq-\delta$ a.e.\ in~$(0,\Tg)$.
Thus, the first part of \eqref{tesi} is proved.
Let us pass to the second one, by assuming that $\Tg<T$.
Then, $\psi(\Tg)=\gamma$, whence $\psi(t)\leq\gamma$ for every 
$t\in[\Tg,T]$\revis{,} by~\eqref{prelim}.
\Edim

At this point, we can continue the proof of Theorem~\ref{Goal}.
We assume that $\rho>\rhostar$ and prepare some material.
Then, we prove the sentences $i)$ and $ii)$ of the statement.

\step
Preliminary remarks

By recalling \accorpa{diffineq}{good}, we apply Lemma~\ref{Diffineq}
to the functions \eqref{defpsiphi} with $\delta=\rho-\rhostar$ and $\gamma=\eps/2$,
and set $\Teps:=\Tg$ if $\dKeps(\thetaz)>\eps/2$.
Next, we observe that \cite[Prop.~2.11, p.~39]{Brezis} implies
\Beq
  \dKeps(\thetaz) \to \dK(\thetaz).
  \label{convdKz}
\Eeq
Now, we recall the strong convergence \eqref{strong} 
and deduce that
\Beq
  \thetaeps(t) \to \theta(t) 
  \quad \hbox{strongly in $H$} \quad \aat .
  \label{aestrong}
\Eeq
Note that we can establish pointwise strong convergence only in the uniformly parabolic case 
(the weak star convergence \eqref{convthetabis} implies strong convergence in $\C0H$, indeed),
while, in general, we just have the almost everywhere strong convergence~\eqref{aestrong}.
We prove that \eqref{aestrong} implies \revis{that}
\Beq
  \dKeps(\thetaeps(t)) \to \dK(\theta(t))
  \quad \aat .
  \label{convdK}
\Eeq
For a fixed $t$ for which the strong convergence \eqref{aestrong} holds,
by setting $\veps:=\thetaeps(t)$ and $v:=\theta(t)$ for brevity,
we have from \pier{\eqref{moryos2}}~\revis{that}
\Beq
  \dK(\theta(t)) - \dKeps(\thetaeps(t))
  = \int_0^{\dK(v)} \bigl( 1 - \min\{s/\eps,1\} \bigr) \, ds
  - \int_{\dK(v)}^{\dK(\veps)} \min\{s/\eps,1\} \, ds.
  \non
\Eeq
The first integral tends to zero by the Lebesgue dominated convergence theorem,
while the absolute value of the second one is estimated by $|\dK(\veps)-\dK(v)|\leq\normaH{\veps-v}$.

\step
First case

Suppose that $\thetaz\in K$.
Then, $\dK(\thetaz)=0$, whence also $\dKeps(\thetaz)=0$.
Thus, $\dKeps(\thetaeps(t))\leq\eps/2$ for every $t\in[0,T]$.
Hence, from~\eqref{convdK}, we infer that $\dK(\theta(t))=0$, i.e., $\theta(t)\in K$, \aat.
In order to extend this to the whole interval~$[0,T]$,
we owe to \eqref{weakcont} and to the fact that the convex set $K$ is weakly closed. 
Therefore, the property $\theta(t)\in K$ proved for a dense subset of $[0,T]$ holds \revis{true} for whole interval.

\step
Second case

Let now suppose that $\thetaz\not\in K$.
Thus $\dK(\thetaz)>0$ and we can assume that $\eps<\dK(\thetaz)$.
Hence, we also have $\dKeps(\thetaz)>\eps/2$
and the time $\Teps\in(0,T]$ is well-defined.
We~set
\Beq
  \Tstar := \liminf_{(\eps,\alpha)\to(0,0)} \Teps \,.
  \label{defTstar}
\Eeq
Thus, $\Tstar\in[0,T]$\revis{,} and we prove that it satisfies the conditions of Theorem~\ref{Goal}.
From the lemma and \eqref{convdKz}, we~have
\Beq
  \Teps \leq \frac {\dKeps(\thetaz)} {\rho-\rhostar} \revis{,}
  \quad \hbox{whence} \quad
  \Tstar \leq \frac {\dK(\thetaz)} {\rho-\rhostar} \revis{,}
\Eeq
i.e., the upper bound of the statement.
We now show that $\Tstar>0$ and argue by contradiction,
i.e., we assume that $\Teps$ tends to zero for a subsequence.
From the inequality $\dKeps(\thetaeps(t))\leq\eps/2$ for $t\geq\Teps$
and \eqref{convdK}, we deduce that $\dK(\theta(t))=0$ \aat.
Hence, we have $\theta(t)\in K$ \aat\ and \eqref{weakcont} yields $\thetaz\in K$, a~contradiction.
Hence, we can consider the non-empty interval $(0,\Tstar)$
and prove~\eqref{decreasing}.
To this end, we take any $\tz\in(0,\Tstar)$.
We can assume that $\Teps>\tz$, so that the first sentence of \eqref{tesi} implies
\Beq
  \frac d{dt} \, \dKeps(\thetaeps(t))) \leq -(\rho-\rhostar) 
  \quad \hbox{for a.a.\ $t\in(0,\tz)$}.
  \label{perdisug}
\Eeq
Now, we fix $\thetabar\in K$ and have 
$\dKeps(\thetaeps(t))\leq\dK(\thetaeps(t))\leq\normaH{\thetaeps(t)-\thetabar}\leq c$ \aat\
thanks to the first \pier{condition in}~\eqref{strong}.
Hence, by accounting for \eqref{convdK} once more, we derive that
$\dKeps(\thetaeps(\cpto))$ converges to $\dK(\theta(\cpto))$ strongly in~$L^2(0,T)$,
whence the convergence in the sense of distributions follows for their derivatives\revis{,}
and the inequality \eqref{perdisug} is conserved in the limit.
As $\tz$ is arbitrary, \eqref{decreasing}~is proved.
Finally, we show~\eqref{goal} by assuming $\Tstar<T$.
Take any $\tz\in(\Tstar,T)$.
We can assume that $\Teps<\tz$ (for a subsequence).
Hence, we have $\dKeps(\thetaeps(t))\leq\eps/2$ for every $t\in[\tz,T]$.
It follows that $\dK(\theta(t))=0$ for a.a.\ $t\in(\tz,T)$\revis{,} by~\eqref{convdK}.
As $\tz$ is arbitrary, we have $\theta(t)\in K$ 
for a.e.\ $t\in(\Tstar,T)$\revis{,}
and we conclude that $\theta(t)\in K$ for every $t\in[\Tstar,T]$\revis{,}
by the weak continuity~\eqref{weakcont}.
This completes the proof.\revis{\qed}


\section{Appendix}
\label{APPENDIX}
\setcounter{equation}{0}

In this section, we prove Lemma~\ref{MorYos}.
We denote by $(\cpto,\cpto)$ the scalar product of~$H$
and write $\norma\cpto$ instead of~$\normaH\cpto$.
We start proving that
the function $\dK$ is \frechet\ differentiable at any point $v\in H\setminus K$
and \revis{that} its gradient is given~by
\Beq
  D\dK(v) = \frac {\QK v} {\normaH{\QK v}} \,.
  \label{derdK}
\Eeq
This easily follows from \cite[Prop.~2.2]{FitzPh}.
Indeed, we immediately deduce \juerg{from this result} that
$2\QK$ is the \frechet\ derivative of the map $v\mapsto\phi(v):=\norma{\QK v}^2$.
Now, we read $\dK$ as the square root of~$\phi$
and assume that $v\in H\setminus K$, i.e., $\phi(v)>0$.
Then, $\phi>0$ in a \nbh\ of $v$, 
whence $\dK$~is \frechet\ differentiable at~$v$
and \eqref{derdK} follows \revis{from} applying the chain rule:
\Beq
  D\dK(v) 
  = \frac 12 \, (\phi(v))^{-1/2} D\phi(v)
  = \frac 12 (\dK(v))^{-1} \, 2\QK v
  = \frac {\QK v} {\norma{\QK v}} \,.
  \non
\Eeq
On the contrary, the proof of the rest of the lemma needs some work.
Assume first $v\in K$.
Then, $v$~is a minimum point for~$\dK$ \pier{and $\dKeps$}, whence $0\in\partial\dK(v)$ and $D\dKeps(v)=0$.
On the other hand, we have $\QK v=0$\juerg{, and thus \pier{\eqref{moryos1} holds} true in 
this case}.
Assume now $\dK(v)>\eps$.
Then, the point
\Beq
  \xi := \frac {\QK(v)} {\dK(v)}
  \non
\Eeq
\pier{satisfies $\norma\xi=1$ and $v-\eps\xi\not\in K$.
Hence, \pier{\eqref{moryos1}} reduces to $D\dKeps(v)=\xi$ 
and thus means that $D\dK(v-\eps\xi)=\xi$ (by \revis{the} definition 
of \revis{the} Yosida regularization).
Therefore}, we prove this fact.
We set $\lambda:=1-\eps/\dK(v)$ and observe that $\lambda>0$ and that
\Beq
  v - \eps\xi - \PK v
  = \QK v - \frac \eps{\dK(v)} \, \QK v
  = \lambda \QK v
  = \lambda (v-\PK v) .
  \label{peryosida}
\Eeq
Then we have\revis{,} for every $z\in K$\revis{,}
\Beq
  (v-\eps\xi-\PK v,z-\PK v)
  = \lambda (v-\PK v,z-\PK v)
  \leq 0.
  \non
\Eeq
As $\PK v\in K$, this shows that $\PK(v-\eps\xi)=\PK(v)$.
By applying \eqref{derdK} to $v-\eps\xi$ and recalling \eqref{peryosida} once more, we obtain
\Beq
  D\dK(v-\eps\xi)
  = \frac {v-\eps\xi-\PK v} {\norma{v-\eps\xi-\PK v}}
  = \frac {\lambda\QK v} {\norma{\lambda\QK v}}
  = \xi .
  \non
\Eeq
Hence, the desired equality is established under the assumption $\dK(v)>\eps$.
As $\dKeps$, $\QK$ and $\dK$ are continuous, 
\pier{\eqref{moryos1}} holds also if $\dK(v)=\eps$,
and we consider the last case, i.e., $0<\dK(v)<\eps$.
We first prove that
\Beq
  \dKeps(v) 
  = \frac 1 {2\eps} \, \pier{(\dK(v))^2} 
  \quad \hbox{if $0<\dK(v)<\eps$}.
  \label{permoryos}
\Eeq
It is well known that the infimum in the definition \eqref{genmoreau} is a minimum.
We first look for a minimum point $z\not\in K$.
Then, in view of~\eqref{derdK}, $z$~has to satisfy
\Beq
  \frac {z-v} \eps + \frac {\QK z} {\norma{\QK z}}  = 0 , \quad
  \hbox{that is} , \quad
  v = z + \eps \, \frac {\QK z} {\norma{\QK z}} \,.
  \non
\Eeq
It easily follows that $\PK v=\PK z$.
Therefore, we have~\revis{that}
\Beq
  \dK(v)
  = \norma{v-\PK v} 
  = \norma{v-\PK z}
  = \Bigl( 1 + \frac \eps {\norma{\QK z}} \Bigr) \norma{z-\PK z}
  = \dK(z) + \eps
  > \eps  \revis{,}
  \non
\Eeq
while we were assuming that $\dK(v)<\eps$.
Therefore, every minimum point $z$ has to belong to~$K$, and we have
\revis{that}
\Beq
  \dKeps(v)
  = \min_{z\in K} \frac 1{2\eps} \, \norma{z-v}^2
  = \frac 1{2\eps} \, \norma{v-\PK v}^2 \revis{,}
  \non
\Eeq
so that \eqref{permoryos} is proved.
Since the set $\graffe{v\in H:\ 0<\dK(v)<\eps}$ is open,
we can differentiate \eqref{permoryos} by applying the chain rule and \eqref{derdK},
and deduce that \pier{\eqref{moryos1}} holds also in this case.
To conclude the proof, we have to derive \pier{\eqref{moryos2}}. \juerg{Now observe that this identity 
trivially holds if $v\in K$ and that $K$ is nonempty. It thus 
suffices} to prove that \juerg{both sides of the identity} have the same \frechet\ gradient.
\juerg{To this end, assume} first that $v\not\in K$.
By differentiating the \rhs\ at $v$ with the chain rule
and applying \eqref{derdK} and \pier{\eqref{moryos1}}, we obtain
\Beq
  \min\bigl\{ \dK(v)/\eps , 1 \bigr\} D\dK(v)
  = \min\bigl\{ \dK(v)/\eps , 1 \bigr\} \, \frac {\QK v} {\dK(v)}
  = \frac {\QK v} {\max\{\eps,\dK(v)\}}
  = D\dKeps(v) .
  \non  
\Eeq
Assume now that $v$ belongs to $K$ and take any $h\in H$ satisfying $\norma h\leq\eps$.
Then, we have $\dK(v+h)\leq\norma h\leq\eps$, and we infer that
\Beq
  0 \leq \int_0^{\dK(v+h)} \min\{s/\eps,1\} \, ds
  = \int_0^{\dK(v+h)} \frac s\eps \, ds
  = \frac 1 {2\eps} \, \pier{(\dK(v+h))^2}
  \leq \frac 1 {2\eps} \, \norma h^2.
  \non
\Eeq
Thus, the \frechet\ gradient of the \rhs\ of  \pier{\eqref{moryos2}} at $v$ is zero.
On the other hand, we also have $D\dKeps(v)=0$ in this case.
This completes the proof.\qed


\section*{Acknowledgments}
PC and GG gratefully acknowledge some financial support from 
\pier{the GNAMPA (Gruppo Nazionale per l'Analisi Matematica, 
la Probabilit\`a e le loro Applicazioni) of INdAM (Istituto 
Nazionale di Alta Matematica) and the IMATI -- C.N.R. Pavia.}


\vspace{3truemm}

\Begin{thebibliography}{10}

\gianni{\bibitem{Barbu}
V. Barbu,,
``Nonlinear Differential Equations of Monotone Types in Banach spaces'',
Springer Monographs in Mathematics. Springer, New York, 2010.}

\pier{%
\bibitem{BCGMR}
V. Barbu, P. Colli, G. Gilardi, G. Marinoschi, E. Rocca,
Sliding mode control for a nonlinear phase-field system,
preprint arXiv:1506.01665~[math.AP] (2015), pp.~1-28.}%

\bibitem{Brezis}
H. Brezis,
``Op\'erateurs maximaux monotones et semi-groupes de contractions
dans les espaces de Hilbert'',
North-Holland Math. Stud.
{\bf 5}.
North-Holland,
Amsterdam,
1973.

\pier{%
\bibitem{CRS11}
M.-B. Cheng, V. Radisavljevic, W.-C. Su,
Sliding mode boundary control of a parabolic PDE system 
with parameter variations and boundary uncertainties, 
{\it Automatica J. IFAC} {\bf 47} (2011), 381-387.
}%

\pier{%
\bibitem{DiB}
E. DiBenedetto, ``Degenerate parabolic equations'',
Universitext. Springer-Verlag, New York, 1993.
}%

\bibitem{FitzPh}
S. Fitzpatrick, R.R. Phelps,
Differentiability of the metric projection in Hilbert space,
Trans. Amer. Math. Soc. {\bf 270} (1982), 483-501.

\pier{%
\bibitem{Levaggi13}
L. Levaggi, 
Existence of sliding motions for nonlinear evolution equations in Banach spaces, 
Discrete Contin. Dyn. Syst. 2013, Dynamical systems, differential equations 
and applications. 9th AIMS Conference. Suppl., 477-487.
\bibitem{LO02}
L. Levaggi,
Infinite dimensional systems' sliding motions,
Eur. J. Control {\bf 8} (2002), 508-516.
}%

\bibitem{Lions}
J.-L.~Lions,
``Quelques m\'ethodes de r\'esolution des probl\`emes
aux limites non lin\'eaires'',
Dunod; Gauthier-Villars, Paris, 1969.

\pier{
\bibitem{PiUs}
A. Pisano, E. Usai,
Sliding mode control: a survey with applications in math,
Math. Comput. Simulation {\bf 81} (2011), 954-979.
}

\pier{%
\bibitem{Roub}
T. Roub{\'{\i}}{\v{c}}ek,  ``Nonlinear partial differential equations with applications''
International Series of Numerical Mathematics {\bf 153}. Birkh\"auser Verlag, Basel, 
2005.}

\pier{%
\bibitem{Show}
R.E. Showalter,
``Monotone Operators in Banach Space and Nonlinear Partial Differential Equations'', 
Mathematical Surveys and Monographs {\bf 49}. American Mathematical Society, 
Providence, RI, 1997.
}%

\pier{%
\bibitem{XLGK13}
H. Xing, D. Li, C. Gao, Y. Kao,
Delay-independent sliding mode control for a class of quasi-linear parabolic distributed parameter systems with time-varying delay, {\it J. Franklin Inst.} {\bf 350} (2013), 397-418.
}%

\pier{%
\bibitem{Zhe}
S. Zheng, ``Nonlinear evolution equations'', 
Chapman \& Hall/CRC Monographs and Surveys in Pure and Applied
Mathematics {\bf 133}. Chapman \& Hall/CRC, Boca Raton, FL, 2004. 		
}

\End{thebibliography}

\End{document}


\step
Lower bound

We recall the definition of $\thetastarn$ given in \eqref{regthetabis}
and prove the uniform lower bound
\Beq
  \thetaeps \geq \thetastarn
  \quad \aeQ
  \label{lower}
\Eeq
which also implies the lower bound stated in \eqref{regthetabis} 
for the function $\theta$ given by~\eqref{convtheta}.
We test \eqref{primaeps} by $\zeta:=-(\thetaeps-\thetastarn)^-$,
where $(\cpto)^-$ denotes the negative part,
and obtain
\Beq
  \frac 12 \iO |\zeta(t)|^2
  + \intQt \keps(\thetaeps) |\nabla\zeta|^2
  + \rho \intQt \sigmaeps \zeta
  = \frac 12 \iO |\zeta(0)|^2
  + \intQt f \zeta .
  \label{perlower}
\Eeq
The \rhs\ is nonpositive since $\thetaz\geq\thetazmin\geq\thetastarn$ and $f\geq0$.
Now, the product $\sigmaeps\zeta$ is zero where $\theta\geq\thetastarn$.
Where $\theta<\thetastarn$, we have $\zeta<0$ and $\thetaeps<\thetastar$,
whence $\QI\thetaeps<0$. 
Thus $\sigmaeps<0$ by \eqref{ptyos} and $\sigmaeps\zeta>0$.
Hence, \eqref{perlower} yields $\zeta=0$ \aeQ, i.e.,~\eqref{lower}.

To this end, we take any pair of points $t,s\in(0,\Tstar)$
for which the strong convergence \eqref{aestrong} holds.
We can assume that $s,t<\Teps$, so that the first sentence of \eqref{tesi} yields
\Beq
  \frac {\dKeps(\thetaeps(t)) - \dKeps(\thetaeps(s))} {t-s} \leq -(\rho-\rhostar) .
  \non
\Eeq
By accounting for \eqref{convdK} once more, we derive the analogue for~$\dK(\theta)$.
Thus
\Beq
  \frac {\dK(\theta(t)) - \dK(\theta(s))} {t-s} \leq -(\rho-\rhostar)
  \quad \hbox{for a.a.\ $t,s\in(0,\Tstar)$}
  \non
\Eeq
and this is equivalent to~\eqref{decreasing}.